\documentclass[12pt]{amsart}
\usepackage{amssymb}
\usepackage{epsf,amsmath}
\usepackage{amssymb,latexsym}

\newtheorem{thm}{Theorem}

\newtheorem{prop}[thm]{Proposition}

\newcommand{\defeq}{\stackrel{\rm def}{=}}
\newcommand{\bd}{\partial}

\def\({\left(}
\def\){\right)}

\def\no={\,{\,|\!\!\!\!\!=\,\,}}
\def\de{\delta}

\def\bb{\bold{b}}

\def\ff{\bold{f}}

\def\hh{\bold{h}}

\def\ep{\varepsilon}
\def\de{\delta}
\def\De{\Delta}

\def\no={\,{\,|\!\!\!\!\!=\,\,}}

\def\({\left(}
\def\){\right)}

\def\beq{\begin{eqnarray}}
\def\eeq{\end{eqnarray}}

\begin{document}
\date{\today}

\title[Chain polynomials of distributive lattices]
{Chain polynomials of distributive lattices are 75 \% unimodal}

\author[Anders Bj\"orner]{Anders Bj\"orner}
\email{bjorner@math.kth.se}
\address{Department of Mathematics\\
         Royal Institute of Technology\\
         S-100~44 Stockholm, Sweden}

\author[Jonathan David Farley]{Jonathan David Farley}

\address{Department of Applied Mathematics\\
        Massachusetts Institute of Technology\\
	Cambridge, Massachusetts 02139\\
        United States of America}

\begin{abstract}
It is shown that the numbers $c_i$ of chains of length $i$
in the proper part $L\setminus\{0,1\}$ of a distributive 
lattice $L$ of length $\ell +2$ satisfy the inequalities
$$c_0<\ldots<c_{\lfloor{\ell /2}\rfloor} \quad\mbox{ and }\quad
c_{\lfloor{3 \ell /4}\rfloor}>\ldots>c_{\ell}.$$
This proves 75 \% of the inequalities implied by the 
Neggers unimodality conjecture.
\end{abstract}

\maketitle

\section{Introduction}

The {\em chain polynomial} of a finite poset $P$ is defined as
$$C(P,t)=\sum_i c_i t^i,
$$
where $c_i$ is the number of chains (totally ordered subsets) in $P$
of length $i$ (i.e., cardinality $i+1$). One of the equivalent forms
of a well-known poset  conjecture due to Neggers \cite{Ne} implies that
the chain polynomial of 
the proper part $L\setminus\{0,1\}$ of a distributive 
lattice $L$ of length $d+1$ is {\em unimodal},
meaning that for some $k$ the coefficients of
$C(L\setminus\{0,1\}  ,t)$ satisfy the inequalities
$$c_0 \le \ldots\le c_k
\ge\ldots\ge c_{d-1}.$$
See \cite{Bre} and \cite{St4} for background, references  and more details
concerning this unimodality conjecture.
Recent progress in a special case (when the poset of join-irreducibles
is graded) appears in  \cite{Bran}, \cite{Far} and \cite{ReWe}.

The purpose of this note is to show that the unimodality
conjecture is 75\% correct, in the sense that violations of
unimodality can occur only for indices (roughly) between $d/2$ and
$3d/4$. More precisely,
we prove the following.

\begin{thm}
The numbers $c_i$ of chains of length $i$
in the proper part of a distributive 
lattice $L$ of length $d+1$ satisfy the inequalities
$$c_0<\ldots<c_{\lfloor{(d-1)/2}\rfloor} \;\mbox{ and }\;
c_{\lfloor{3(d-1)/4}\rfloor}>\ldots>c_{d-1}.$$
\end{thm}

The proof consists in observing that the order complex of
$L\setminus\{0,1\}$ is a nicely behaved ball, and then gathering and 
combining some known facts from $f$-vector theory. The pieces of the
argument are stated as Propositions 2, 3, 4 and 5. Of these, only
Proposition 3 seems to be new.

\section{Some $f$-vector inequalities}

For standard notions concerning simplicial
complexes we refer to the literature, see e.g. the books \cite{Br, Z}.

Let $\De$ be a $(d-1)$-dimensional simplicial complex, and let $f_i$ be the
number of $i$-dimensional faces of $\De$. The sequence $(f_0,\dots,f_{d-1})$
is called the {\em $f$-vector} of $\De$. We put $f_{-1}=1$.
The {\em $h$-vector} $(h_0,\dots,h_d)$ of $\De$ is defined by the equation
\begin{equation}\label{fh-rel}
\sum_{i=0}^{d} f_{i-1} x^{d-i} = \sum_{i=0}^{d} h_i (x+1)^{d-i}. 
\end{equation}
\bigskip

In the following two results we 
assume that  $\(f_0,f_1,\ldots,f_{d-1}\)$ is the
$f$-vector of a $(d-1)$-dimensional simplicial complex $\De$,
and that $f_0 >d$.
From now on, let $d\ge 3$ and $\de\defeq \lfloor\frac{d}2\rfloor$,
$\ep\defeq \lfloor\frac{d-1}2\rfloor$.

\begin{prop} Suppose that 
$h_i \ge 0,$ for all $0\le i \le d$.
Then
$$ f_i < f_j, \mbox{   for all $i<j$ such that $i+j\le d-2.$}$$
In particular, $f_0<f_1<\ldots<f_{\ep}$.
\end{prop}

\begin{proof}
This implication is well known. See e.g. \cite[Proposition 7.2.5 (i)]{Bj1}.
\end{proof}

\begin{prop} Suppose that $h_i \ge h_{d-i} \ge 0, \quad \mbox{for all $0\le i \le \de$.}$
Then
$$f_{\lfloor{3(d-1)/4}\rfloor}>
\ldots>f_{d-2}>f_{d-1}.$$
\end{prop}

\medskip

\begin{proof}
By (\ref{fh-rel}),
the $f$-vector $\ff=\(f_0,f_1,\ldots,f_{d-1}\)$ 
and the $h$-vector $\hh=\(h_0,h_1,\ldots,h_d\)$ satisfy
\begin{equation}\label{fh2}
f_k=\sum^d_{i=0}h_i\binom{d-i}{d-1-k}\qquad , \qquad k=-1,\ldots,d-1.
\end{equation}
Define integer vectors $\bb^i$ as follows:
$$
\bb^i=\(b^i_0,b^i_1,\ldots,b^i_{d-1}\) \qquad, \qquad\text{where }
b^i_k=\binom{i}{d-1-k}.
$$
Then, by  (\ref{fh2}), 
$\ff=\sum^d_{i=0}h_i\bb^{d-i}$,
which we rewrite
\begin{equation}\label{fh3}
\ff=\sum^{\ep}_{i=0}(h_i -h_{d-i})\bb^{d-i} +  \sum^{\de}_{i=0}h_{d-i}{\tilde{\bb}}^i,
\end{equation}
where
$$
{\tilde{\bb}}^i\defeq\left\{
\begin{array}{ll}
\bb^i+\bb^{d-i}&\qquad,\qquad\text{if } 2i\no=d\\
\bb^{d/2}&\qquad,\qquad\text{if } 2i=d. 
\end{array}\right.
$$

Let us say that a unimodal sequence
$$
a_0\le a_1\le\ldots\le a_k\ge a_{k-1}\ge\ldots\ge a_n
$$
{\it peaks\/} at $k$ 
(note that this does not necessarily determine $k$ uniquely).

It is shown in \cite[Proof of Thm. 5, p. 50]{Bj2} that
the vector ${\tilde{\bb}}^i$ is uni\-modal and peaks at
$d-1-\lfloor\frac{(d-i)}2\rfloor$.
The vector $\bb^{d-i}$ is a segment of a row in Pascal's triangle,
so it is easy to see that it is uni\-modal and, in fact,  also peaks at
$d-1-\lfloor\frac{(d-i)}2\rfloor$. One easily checks that
$$
d-1-\lfloor\frac{(d-i)}2\rfloor =
\left\{
\begin{array}{ll}
\lfloor\frac{d}2\rfloor+\lfloor{\frac{i}2}\rfloor-1&\qquad,\qquad\text{if
$d$ and $i$ are even}\\
\lfloor{\frac{d}2}\rfloor+\lfloor{\frac{i}2}\rfloor&\qquad,\qquad
\text{otherwise.}
\end{array}\right.
$$
Hence, both the vectors $\bb^{d-i}$ ($0 \le i\le \ep$) and the vectors
 $\tilde{\bb}^i$ ($0 \le i\le \de$) are unimodal and peak between $\de$
and $\de+\lfloor\de/2\rfloor$. 

By equation (\ref{fh3}), $\ff$ is a nonnegative linear
combination of the vectors $\bb^{d-i}$ and $\tilde{\bb}^i$.
It follows from the previous paragraph that
the inequalities hold for each of these vectors separately, 
strictly for $\bb^{d}$, and 
non-strictly otherwise. For the computation of the index
$\lfloor{3(d-1)/4}\rfloor$, see again
 \cite[pp. 50--51]{Bj2}.
 Hence, if $h_d=0$ the result follows. 
The case when $h_d=1$ requires a small extra argument to see that the
inequalities are in fact strict. For this case one can proceed as
in \cite[Proof of Thm. 5]{Bj2}.
\end{proof}

\section{On the $h$-vectors of balls}

We say that a simplicial complex is a {\em polytopal $(d-1)$-sphere} if it
is combinatorially isomorphic to the boundary complex of some convex $d$-polytope.
See Ziegler \cite{Z} for notions relating to 
polytopes and convex geometry.

\smallskip

We now review some definitions and results from the general theory of face 
numbers. For more about this topic, see e.g.   
\cite{Z} or the survey \cite{BB}.

It follows from (\ref{fh-rel}) that $h_0=1$, $h_1=f_0 -d$, and $h_d=(-1)^{d-1} \tilde{\chi}(\De)$, where
$\tilde{\chi}(\De)$ is the reduced Euler characteristic of $\De$. In
particular,
$$h_d=\left\{ \begin{array}{ll}
1, & \mbox{if $\De$ is a sphere,} \\
0, & \mbox{if $\De$ is a ball,}
\end{array}
\right.
$$
where the conditions are shorthand for saying that $\De$'s geometric
realization is homeomorphic to a sphere, resp. a ball.

The following are the {\em Dehn-Sommerville relations}:
\begin{equation}\label{DS}
\mbox{{\em If $\De$ is a sphere then \;
$h_i = h_{d-i}$, \; for all $0\le i \le d$.}}
\end{equation}
Hence, for spheres all $f$-vector information is encoded in the
shorter {\em $g$-vector}
$g=(g_0,\dots,g_{\lfloor \frac{d}{2} \rfloor})$, defined by
$g_i =h_i - h_{i-1}$. 
The relevance of the $g$-vector for this paper
is the following result, due to Stanley \cite{St1}:
\begin{equation}\label{g-thm}
\mbox{{\em If $\De$ is a polytopal sphere, then \:
$g_i\geq 0$ for all $i\geq 0$.}}
\end{equation}

If $\De$ is a $(d-1)$-ball, its boundary complex $\bd\De$ is a $(d-2)$-sphere.
Furthermore, $\bd\De$'s  $f$-vector is determined by that of $\De$, as shown by the
following consequence of the Dehn-Sommerville relations, due to
McMullen and Walkup \cite{MW}, see also \cite[Coroll. 3.9]{BL2}:
\begin{equation}\label{ball-thm}
\mbox{{\em If $\De$ is a ball with boundary $\bd \De$, then \;
$h^{\De}_{i} - h^{\De}_{d-i} = g^{\bd\De}_{i}$.}}
\end{equation}

Say that a $(d-1)$-ball $\De$ admits a {\em polytopal embedding} if
$\De$ is isomorphic to a subcomplex of the boundary complex
of some simplicial $d$-polytope. The following was shown by Kalai \cite[\S 8]{Kal}
and Stanley \cite[Coroll. 2.4]{St3}.
\begin{equation}\label{g-thm2}
\mbox{{\em If $\De$ admits a polytopal embedding, then \:
$g^{\bd\De}_i\geq 0$ for all $i\geq 0$.}}
\end{equation}

Combining (\ref{g-thm}), (\ref{ball-thm}) and (\ref{g-thm2}), we deduce 
 the following result.

\begin{prop}
If $\De$ is a $(d-1)$-ball, such that either the 
boundary sphere $\bd\De$ is polytopal
or $\De$ admits a polytopal embedding,
then $$h_i \ge h_{d-i} \ge 0, \; \mbox {for all} \; 0\le i \le \de.$$
\hfill$\Box$
\end{prop}

\section{Proof of Theorem 1}

We refer to \cite[Ch. 3]{St2} for basic facts and notation concerning
distributive lattices.

Let $L$ be a distributive lattice of length $d+1$, and let 
$\De_L =\De(L\setminus \{0,1\})$
be the order complex of its proper part. Thus, $\De_L$ is a 
pure simplicial complex of dimension $d-1$.

\newpage

\begin{prop} 
Suppose that $L$ is 
not Boolean. Then
the complex $\De_L$ is a $(d-1)$-ball satisfying
\begin{enumerate}
\item[(i)] $\De_L$ admits a polytopal embedding,
\item[(ii)] $\bd\De_L$ is polytopal.
\end{enumerate}
\end{prop}

\begin{proof}
By Birkhoff's representation theorem (see \cite[Ch. 3]{St2}) we have
that $L=J(P)$, where $J(P)$ is the family of order ideals of some poset $P$
ordered by inclusion.
Let $B$ denote the Boolean lattice of {\em all} subsets of $P$. Then
$\De_B =\De(B\setminus \{0,1\})$
is a polytope boundary (the barycentric subdivision of
the boundary of a $d$-simplex). Furthermore, $\De_L$ is embedded
in $\De_B$ as a full-dimensional subcomplex. Finally,
$\De_L$ is a shellable ball \cite{Bj3, Pro}. Thus,
part (i) is proved.

Part (ii) requires a small convexity argument. Alternatively, it follows
from Provan's result \cite{Pro} that $\De_L$ can be obtained from a simplex
via repeated stellar subdivisions. Since this part is not needed for
the main result of this paper, details of the proof are left out.
\end{proof}

We now have all the pieces needed to prove Theorem 1.
We may assume that $L$ is not Boolean, 
since in that case $\De_L$ is a sphere and Theorem 1 
is a special case of \cite[Thm. 5]{Bj2}. Then, 
by Propositions 4 and 5 we have that
$$h_i \ge h_{d-i} \ge 0, \; \mbox {for all} \; 0\le i \le \de.$$
Furthermore, by Propositions 2 and 3 it follows that
the $f$-vector of $\De_L$ satisfies
$$f_0<\ldots<f_{\lfloor{(d-1)/2}\rfloor} \;\mbox{ and }\;
f_{\lfloor{3(d-1)/4}\rfloor}>\ldots>f_{d-1}.$$
Since $f_i=c_i$ for all $i$, the proof of Theorem 1 
is complete.

\end{document}